\documentclass[twocolumn]{autart}    % Enable this line and disable the
 \usepackage{mathrsfs}
\usepackage{amsfonts}
\usepackage{stmaryrd}
\usepackage{amsmath, amssymb}
\usepackage{cases}                                    % document whose style resembles the
\usepackage{colortbl}
\usepackage{graphicx}                                     % preceding line to obtain a two-column
\usepackage{graphicx}          % Include this line if your
\begin{document}

\begin{frontmatter}

\title{Stabilization Control for It$\hat{o}$ Stochastic System with Indefinite State and Control Weight Costs\thanksref{footnoteinfo}} % Title, preferably not more

\thanks[footnoteinfo]{This work is supported by the National Natural Science Foundation of
China under Grants 61120106011, 61573221, 61633014. Corresponding author Huanshui Zhang.}

\author{Hongdan Li}\ead{lhd200908@163.com},    % Add the
\author{Qingyuan Qi}\ead{qiqy123@163.com},               % e-mail address
\author{Huanshui Zhang}\ead{hszhang@sdu.edu.cn}  % (ead) as shown

\address{School of Control Science and Engineering,
Shandong University, Jinan, Shandong, P.R.China 250061.}  % Please supply

\begin{keyword}                           % Five to ten keywords,
Indefinite stochastic linear quadratic control; generalized differential Riccati equation; convergence; stabilization.            % chosen from the IFAC
\end{keyword}                             % keyword list or with the
                                          % help of the Automatica
                                          % keyword wizard

\begin{abstract}
In standard linear quadratic (LQ) control,  the first step in investigating infinite-horizon optimal control is to derive the stabilization condition with the optimal LQ controller.  This paper focuses on the stabilization of an It$\hat{o}$ stochastic system with indefinite control and state weighting matrices in the cost functional. A generalized algebraic Riccati equation (GARE) is obtained via the convergence of the generalized differential Riccati equation (GDRE) in the finite-horizon case. More importantly, the necessary and sufficient stabilization conditions for indefinite stochastic control are obtained.\\
One of the key techniques is that the solution of the GARE is decomposed into a positive semi-definite matrix that satisfies the singular algebraic Riccati equation (SARE) and a constant matrix that is an element of the set satisfying certain  linear matrix inequality conditions. Using the equivalence between the GARE and SARE, we reduce the stabilization of the general indefinite case to that of the definite case, in which the stabilization is studied using a Lyapunov functional defined by the optimal cost functional  subject to the SARE.
\end{abstract}

\end{frontmatter}

\section{Introduction}
The linear quadratic (LQ) control pioneered by Kalman \cite{1Kalman:60} is a classical yet important problem in both theory and engineering applications. In 1968, Wonham \cite{2Wonham:68} investigated stochastic LQ problems, and this topic has since been studied by many researchers \cite{6Zhang:2017}, \cite{5Davis:77}. Most results were obtained under the common assumption that the state weighting matrices are non-negative definite and the control weighting matrices are positive definite.\\
However, in \cite{7Chen:98}, this common assumption was changed, i.e., stochastic LQ problems with well-posed indefinite control weighting matrices were considered. This phenomenon is known as an indefinite stochastic problem, and it has a deeply uncertain nature; for a more detailed discussion and many examples, see \cite{7Chen:98}. As they have a wide range of applications, from portfolio selection to pollution control, an increasing number of researchers have studied indefinite stochastic problems.  For example, \cite{8Rami:2001} relaxed the positive definiteness constraint in \cite{7Chen:98} and solved the indefinite LQ problem by introducing a generalized differential Riccati equation (GDRE). Under the assumption that the system is stabilizable, \cite{9Wu:2002} showed that the solvability of indefinite stochastic LQ problems in the infinite horizon is equivalent to the existence of a static stabilizing solution to the generalized algebraic Riccati equation (GARE). With regard to discrete cases, \cite{10Ferrante:2015} derived an equivalent condition for the discrete-time indefinite optimal control problem in finite horizon. The  LQ problem for discrete time-invariant systems  with arbitrary terminal weight was explored in \cite{11Bilardi:2007}. \cite{12Ni:2017} discussed  an indefinite stochastic LQ problem with state transmission delay and multiplicative noise. Moreover, because indefinite stochastic LQ problems can be understood as the dual versions of robust filtering problems, some indefinite LQ results have been applied to filtering problems   \cite{13Zorzi:2017}, \cite{14Zorzi:2017}.\\
The aforementioned papers mainly studied the LQ optimal control problem, whereas relatively few studies have concentrated on the stabilization problem for an indefinite stochastic problem. Nevertheless, the optimal control problem in the infinite horizon case is a worthy topic of research when the system is stabilizable. For instance, \cite{15Zhang:2018} and \cite{16Qi:2017} investigated stabilization control for linear discrete- and continuous-time mean-field systems, respectively; in \cite{17Rami:2000}, \cite{9Wu:2002} and \cite{18Li:2003}, the basic assumption that the system is mean-square stabilizable was imposed throughout.  As stated in \cite{8Rami:2001}, indefinite stabilization  in the infinite-horizon case is a crucial issue.  In this paper, the indefinite stabilization problem for an It$\hat{o}$ stochastic system is investigated. \\
The main contributions of this paper are as follows.
First, we consider the convergence of the GDRE involving a matrix pseudo-inverse and two additional equality/inequality constraints, which is a weaker requirement than in previous work \cite{19Rami:2001}. In fact, it is natural to relax the matrix invertibility constraint, because one cannot generally know in advance whether a singularity will occur.  Second, in \cite{19Rami:2001},  the asymptotic behavior of the GDRE was only investigated under a strict positive-definite constraint, and the corresponding stabilization results for the system were not given. In contrast, this paper discusses mean-square stabilization for the It$\hat{o}$ stochastic system. The key technique is to decompose the solution of the GARE into a positive semi-definite matrix that satisfies the singular algebraic Riccati equation (SARE) and a constant matrix that is an element of the set satisfying certain linear matrix inequality (LMI) conditions. In view of the equivalence between the GARE and the SARE,  the stabilization of the general indefinite case can be reduced to that of the definite case, where the stabilization is studied using a Lyapunov functional defined with the optimal cost functional  subject to the SARE.\\
The remainder of this paper is organized as follows.  Some useful preliminary results are given in Section 2. In Section 3, we discuss the convergence of the GDRE and the mean-square stabilization problem. A two-dimensional numerical example is presented in Section 4, and a summary  is provided in Section 5.

\section{Preliminaries}
Let $\mathcal{R}^{n}$ be the family of $n$-dimensional vectors; $M'$ is the transpose of $M$ and $M^{\dagger}$ is the Moore-Penrose pseudo-inverse of $M$;
$M > 0 (\geq0)$ denotes a symmetric matrix  that is strictly positive-definite (positive semi-definite); $\mathbf{Ker}(M)$ is the kernel of a matrix $M$. $(\Omega,\mathcal{F},P,\mathcal{F}_{t} | t\geq0 )$ is
a complete stochastic basis such that $\mathcal{F}_{ 0}$ contains all $P$-null
elements of $\mathcal{F}$, and the filtration is generated by the standard
Brownian motion $\{w(t)\}_{t\geq0}$.
\begin{eqnarray*}
L_{\mathcal{F}}^{2}(0,T;\mathcal{R}^{m})&=&
\Big\{\varphi(t)_{t\in[0,T]}\ \text {is an} \ \mathcal{F}_{t} \ adapted \ stochastic \\
&& process \ s.t.\ \ E\int_{0}^{T}\|\varphi(t)\|^{2}dt<\infty\Big\}.
\end{eqnarray*}\\
We consider the following linear It$\hat{o}$ stochastic system:
\begin{eqnarray}     %方程组开始
\left\{                        %方程组的左边包括大括号\{
\begin{array}{lll}       %设定列阵的格式：{lll}是三个L，表示三列的对齐方式为Left对齐
dx(t)=[Ax(t)+Bu(t)]dt+[Cx(t)+Du(t)]dw(t), \\  %$――分隔列的标记，\\――表示换行
 x(0)=x_{0},%$同上

\end{array}              %方程列阵的结束
\right. \label{f01}                      % 方程组的右边无符号，利用“.“来标示
\end{eqnarray}        %方程组结束
where $x(t) \in \mathcal{R}^{n} $ is the state; $u(\cdot)$, the
admissible control, is any element in $\mathcal{U}_{ad}\equiv L_{\mathcal{F}}^{2}(0,T;\mathcal{R}^{m})$;
and $w(t)$ is the one-dimensional standard Brownian motion. $x_{0} \in \mathcal{R}^{n}$ is the initial value, and $A,B,C,D$ are constant matrices with
compatible dimensions.\\
The associated quadratic cost functional with an infinite horizon is
\begin{eqnarray}
J(u(\cdot))=\frac{1}{2}E\int_{0}^{\infty}\left[x'(t)Q x(t)+u'(t)R u(t)\right]dt,\label{f02}
\end{eqnarray}
where $Q, R$  are symmetric matrices with compatible dimensions.\\
{\bf Problem 1} \ \ Find the optimal control $u(t)=Kx(t)$ with constant
matrix gain $K$ that stabilizes (\ref{f01}) while minimizing (\ref{f02}).\\
To solve Problem 1, we define the finite-horizon control as follows:\\
Consider system (\ref{f01}) with the following cost functional in the finite horizon:
\begin{eqnarray}
J(u(\cdot))&=&\frac{1}{2}E\int_{0}^{T}\left[x'(t)Q x(t)+u(t)'R u(t)\right]dt\nonumber\\
&&+x'(T)P(T)x(T), \label{f03}
\end{eqnarray}
where $Q, R, P(T)$  are symmetric matrices with compatible dimensions.\\
The following result shows the necessary and sufficient condition for the optimal control problem in the finite horizon.\\
{\bf Lemma 1 \cite{8Rami:2001}} \ \ The following conditions are equivalent.\\
(1) \ \ There exists an optimal control $u(t)\in \mathcal{U}_{ad}$ that minimizes (\ref{f03}) subject to (\ref{f01});\\
(2) \ \ The following GDRE admits a solution for $t\in [0,T]$:
 \begin{eqnarray}
\left\{
\begin{array}{lll}
\dot{P}(t)+A'P(t)+P(t)A+C'P(t)C+Q\\
-(P(t)B+C'P(t)D)(R+D'P(t)D)^{\dagger} \\
\times(B'P(t)+D'P(t)C)=0, \\
(R+D'P(t)D)(R+D'P(t)D)^{\dagger} \\
\times(B'P(t)+D'P(t)C)
=B'P(t)+D'P(t)C,\\
R+D'P(t)D\geq 0,
\end{array}
\right.\label{f04}
\end{eqnarray}
with terminal values $P(T)$.  In this case, the optimal controller is $u(t)=K(t)x(t)$, where $K(t)=-(R+D'P(t)D)^{\dagger}(B'P(t)+D'P(t)C)$.
 Moreover, the value functional is
 $\inf_{u(\cdot)}J(u(\cdot))=x_{0}'P(0)x_{0}$.\\
 Based on the above result in the finite time horizon, we will show the main results in the next section.
%Before presenting the main results of the paper,  we first list some lemmas as in the following.\\

\section{Main results}
In this section, we  discuss the convergence of the GDRE (\ref{f04}) and the mean-square stabilization problem of system (\ref{f01}).\\
Define the GARE as follows:
\begin{eqnarray}
\left\{
\begin{array}{lll}
A'P+PA+C'PC+Q-(PB+C'PD)\\
\times(R+D'PD)^{\dagger}(B'P+D'PC)=0,\\
(R+D'PD)(R+D'PD)^{\dagger}(B'P+D'PC)\\
=B'P+D'PC, \\
R+D'PD\geq 0.
\end{array}
\right.\label{f06}
\end{eqnarray}\\
{\bf Remark 1 }  Note that the defined Riccati equation (\ref{f06}) is different from that  in \cite{19Rami:2001}, where $R+D'PD$ must be positive definite.\\
Consider the following  set $\mathcal{P}$, which involves the LMI
\begin{eqnarray*}
\mathcal{P}\triangleq
\Bigg\{
\begin{smallmatrix}
\hat{P}=\hat{P}'\Bigg|
&&\left[
  \begin{smallmatrix}
    A'\hat{P}+\hat{P}A+C'\hat{P}C+Q & \hat{P}B+C'\hat{P}D\\
    B'\hat{P}+D'\hat{P}C & R+D'\hat{P}D\\  %第二行元素
  \end{smallmatrix}
\right]\geq 0,\\
&&\mathbf{Ker}(R+D'\hat{P}D)\subseteq (\mathbf{Ker} B \cap \mathbf{Ker} D)
\end{smallmatrix}\Bigg\}
\end{eqnarray*}\\
{\bf Definition 1} \ \  A solution to the GARE (\ref{f06}) is called a maximal solution, denoted by $P_{max}$, if
\begin{eqnarray*}
P_{max}\geq \hat{P}, \quad \forall \hat{P}\in\mathcal{P}.
\end{eqnarray*}
For convenience, we denote $P(t)$ as $P(t,T)$ with the terminal time $T$ in the GDRE (\ref{f04}) and the terminal value $P(T,T)=\hat{P}$ with $\hat{P}\in\mathcal{P}$.\\
{\bf Theorem 1} \ \ Assume $\mathcal{P}\neq \emptyset$. If system (\ref{f01}) is mean-square stabilizable, then $P(t,T)$ is convergent as $t\rightarrow-\infty$ and the limit of $P(t,T)$ denoted by $\bar{P}$ is the maximal solution to the GARE (\ref{f06}).\\
{ \emph{Proof}:} \ \ Let $\hat{P}\in\mathcal{P}$ and define
\begin{eqnarray}
\left\{
\begin{array}{lll}
Q_{\hat{P}}=A'\hat{P}+\hat{P}A+C'\hat{P}C+Q,\\
L_{\hat{P}}=\hat{P}B+C'\hat{P}D,\\
R_{\hat{P}}=R+D'\hat{P}D.
\end{array}
\right.\label{f07}
\end{eqnarray}
Thus, from $\left[
  \begin{array}{cc}
    Q_{\hat{P}}& L_{\hat{P}}\\
    L_{\hat{P}}' & R_{\hat{P}}\\  %第二行元素
  \end{array}
\right]\geq 0$ and Schur's Lemma, we have  $$R_{\hat{P}}\geq 0,\quad  Q_{\hat{P}}-L_{\hat{P}}R_{\hat{P}}^{\dagger}L_{\hat{P}}'\geq 0, \quad L_{\hat{P}}(I-R_{\hat{P}}R_{\hat{P}}^{\dagger})=0.$$\\
On this basis, i.e., $\left[
  \begin{array}{cc}
    Q_{\hat{P}}& L_{\hat{P}}\\
    L_{\hat{P}}' & R_{\hat{P}}\\  %第二行元素
  \end{array}
\right]\geq 0$ and $R_{\hat{P}}\geq 0$,  we consider the following singular differential Riccati equation (SDRE) with terminal values $Z_{\hat{P}}(T, T)=0$:
\begin{eqnarray}
&&\dot{Z}_{\hat{P}}(t,T)+A'Z_{\hat{P}}(t,T)+Z_{\hat{P}}(t,T)A+C'Z_{\hat{P}}(t,T)C
+Q_{\hat{P}}\nonumber\\
&&-(Z_{\hat{P}}(t,T)B+C'Z_{\hat{P}}(t,T)D+L_{\hat{P}}) (R_{\hat{P}}+D'Z_{\hat{P}}(t,T)D)^{\dagger}\nonumber\\
&&\times(B'Z_{\hat{P}}(t,T)+D'Z_{\hat{P}}(t,T)C+L_{\hat{P}}')=0,\label{f08}
\end{eqnarray}
with the following regular condition
\begin{eqnarray}
&&[I-(R_{\hat{P}}+D'Z_{\hat{P}}(t,T)D)(R_{\hat{P}}+D'Z_{\hat{P}}(t,T)D)^{\dagger}]\nonumber\\
&&\times(B'Z_{\hat{P}}(t,T)+D'Z_{\hat{P}}(t,T)C+L_{\hat{P}}')=0.\label{f09}
\end{eqnarray}
Using the classical ordinary differential equation
theory, it is easy to see that there exists a positive semi-definite solution  to (\ref{f08}). Denote by $Z_{\hat{P}}(\cdot, T)$ the solution of (\ref{f08}) with $Z_{\hat{P}}(T, T)=0$.\\
Next, we mainly investigate whther $Z_{\hat{P}}(\cdot, T)$ also satisfies (\ref{f09}). From the singular value decomposition of $R_{\hat{P}}+D'Z_{\hat{P}}D$ (see Theorem 2.6.3 in \cite{23Horn:90}), there exists an orthogonal matrix $U$ satisfying
\begin{eqnarray}
(R_{\hat{P}}+D'Z_{\hat{P}}D)^{\dagger}=U\left[
  \begin{array}{cc}
    \Lambda^{-1} & \textbf{0} \\
    \textbf{0}& \textbf{0}
  \end{array}
\right]U',\label{f10}
\end{eqnarray}
where $\Lambda$ is an invertible matrix that has the same dimensions as the rank of $R_{\hat{P}}+D'Z_{\hat{P}}D$.
Now, decompose $U$  as
$\left[
  \begin{array}{cc}
U_{1} & U_{2}
  \end{array}
\right]$. The orthogonality of $U$ implies that the columns of the matrix $U_{2}$ form a basis of $\mathbf{Ker}(R_{\hat{P}}+D'Z_{\hat{P}}D)$. The positive semi-definite nature of matrices $R_{\hat{P}}, Z_{\hat{P}}$ implies that $\mathbf{Ker}(R_{\hat{P}}+D'Z_{\hat{P}}D)\subseteq \mathbf{Ker}(R_{\hat{P}})$. In view of $\mathbf{Ker}(R_{\hat{P}})\subseteq \mathbf{Ker}(B)\cap \mathbf{Ker}(D)$,  we have that
\begin{eqnarray}
&&(B'Z_{\hat{P}}+D'Z_{\hat{P}}C+L_{\hat{P}})'[I-(R_{\hat{P}}+D'Z_{\hat{P}}D)\nonumber\\
&&\times(R_{\hat{P}}+D'Z_{\hat{P}}D)^{\dagger}]\nonumber\\
&=&(B'Z_{\hat{P}}+D'Z_{\hat{P}}C+L_{\hat{P}})'U_{2}U_{2}'\nonumber\\
&=&[(Z_{\hat{P}}+\hat{P})B+C'(Z_{\hat{P}}+\hat{P})D]U_{2}U_{2}'\nonumber\\
&=&\textbf{0}.\label{f010}
\end{eqnarray}
 Therefore, $Z_{\hat{P}}(\cdot, T)$ is a solution of the SDRE (\ref{f08})--(\ref{f09}).\\
Considering the following cost functional corresponding to  system (\ref{f01}):
\begin{eqnarray}
\tilde{J}_{T}(u(\cdot))=E\int_{0}^{T}\left[
  \begin{array}{cc}
x(t)\\
u(t)
  \end{array}
\right]'\left[
  \begin{array}{cc}
Q_{\hat{P}} &L_{\hat{P}}\\
L_{\hat{P}}' & R_{\hat{P}}
  \end{array}
\right]\left[
  \begin{array}{cc}
x(t)\\
u(t)
  \end{array}
\right]dt,\label{f11}
\end{eqnarray}
from $\mathcal{P}\neq \emptyset$, we know that $\tilde{J}_{T}(u(\cdot))\geq 0$.\\
Applying Lemma 1,  we find that
\begin{eqnarray}
\tilde{V}_{T}\triangleq\inf\limits_{u(\cdot)}\tilde{J}_{T}(u(\cdot))
=x_{0}'Z_{\hat{P}}(0,T)x_{0}.\label{f12}
\end{eqnarray}
Note the time-invariance of (\ref{f08}) with respect to $T$, with terminal time $T$ and terminal value $Z_{\hat{P}}(T,T)=0$. Thus, for any $t_{1}<t_{2}\leq T$ and for all $x_{0}\neq 0$, we have
\begin{eqnarray*}
&\tilde{V}_{T-t_{1}}=x_{0}'Z_{\hat{P}}(0,T-t_{1})x_{0}=x_{0}'Z_{\hat{P}}(t_{1}, T)x_{0}\\
\geq&\tilde{V}_{T-t_{2}}=x_{0}'Z_{\hat{P}}(0,T-t_{2})x_{0}=x_{0}'Z_{\hat{P}}(t_{2}, T)x_{0}.
\end{eqnarray*}
Because  $x_{0}$ is arbitrary, we have that $Z_{\hat{P}}(t_{1}, T)\geq Z_{\hat{P}}(t_{2}, T)$. Similarly, when $t\leq T_{1}<T_{2}$,
\begin{eqnarray*}
&\tilde{V}_{T_{1}-t}=x_{0}'Z_{\hat{P}}(0,T_{1}-t)x_{0}=x_{0}'Z_{\hat{P}}(t, T_{1})x_{0}\\
\leq&\tilde{V}_{T_{2}-t}=x_{0}'Z_{\hat{P}}(0,T_{2}-t)x_{0}=x_{0}'Z_{\hat{P}}(t, T_{2})x_{0}.
\end{eqnarray*}
 That is, $Z_{\hat{P}}(t, T_{1})\leq Z_{\hat{P}}(t, T_{2})$. Therefore, $Z_{\hat{P}}(t, T)$ is monotonically
increasing with respect to $T$ and is monotonically decreasing
with respect to $t$.\\
Next, we show that $Z_{\hat{P}}(t, T)$ is bounded. The mean-square stabilizability of system (\ref{f01}) with the controller $u(t)=Kx(t)$ yields
\begin{eqnarray*}
x_{0}'Z_{\hat{P}}(0,T)x_{0}\leq \tilde{J}_{T}&\leq& E\int_{0}^{\infty}\left[
  \begin{array}{cc}
x(t)\\
u(t)
  \end{array}
\right]'\left[
  \begin{array}{cc}
Q_{\hat{P}} &L_{\hat{P}}\\
L_{\hat{P}}' & R_{\hat{P}}
  \end{array}
\right]\left[
  \begin{array}{cc}
x(t)\\
u(t)
  \end{array}
\right]dt\\
&=&E\int_{0}^{\infty}x(t)'\left[
  \begin{array}{cc}
I\\
K
  \end{array}
\right]'\left[
  \begin{array}{cc}
Q_{\hat{P}} &L_{\hat{P}}\\
L_{\hat{P}}' & R_{\hat{P}}
  \end{array}
\right]\left[
  \begin{array}{cc}
I\\
K
  \end{array}
\right]\\
&&\cdot x(t)dt\\
&\leq& \lambda_{max}E\int_{0}^{\infty}x'(t)x(t)dt<+\infty,
\end{eqnarray*}
 where $\lambda_{max}$ denotes the maximum eigenvalue of a matrix $\left\{\left[
  \begin{array}{cc}
I\\
K
  \end{array}
\right]'\left[
  \begin{array}{cc}
Q_{\hat{P}} &L_{\hat{P}}\\
L_{\hat{P}}' & R_{\hat{P}}
  \end{array}
\right]\left[
  \begin{array}{cc}
I\\
K
  \end{array}
\right]\right\}$. With an arbitrary  $x_{0}$, we find that $Z_{\hat{P}}(t,T)$ is bounded. Hence, we have that
$Z_{\hat{P}}(t,T)$ is convergent, i.e.,
\begin{eqnarray*}
\lim\limits_{t\rightarrow-\infty}Z_{\hat{P}}(t,T)
=\lim\limits_{T\rightarrow+\infty}Z_{\hat{P}}(0,T-t)\triangleq \bar{Z}_{\hat{P}},
\end{eqnarray*}
where $\bar{Z}_{\hat{P}}$ is a positive semi-definite constant matrix that is independent of $t$.\\
Let $t\rightarrow -\infty$ in the SDRE (\ref{f08})--(\ref{f09}). We have that $\bar{Z}_{\hat{P}}$ is a solution of the following SARE:
\begin{eqnarray}
\left\{
\begin{array}{lll}
A'\bar{Z}_{_{\hat{P}}}+\bar{Z}_{_{\hat{P}}}A+C'\bar{Z}_{_{\hat{P}}}C+Q_{_{\hat{P}}}\\
-(\bar{Z}_{_{\hat{P}}}B+C'\bar{Z}_{_{\hat{P}}}D+L_{_{\hat{P}}})(R_{_{\hat{P}}}+D'\bar{Z}_{_{\hat{P}}}D)^{\dagger}\\
\times(B'\bar{Z}_{_{\hat{P}}}+D'\bar{Z}_{_{\hat{P}}}C+L_{_{\hat{P}}}')=0,\\
 (R_{_{\hat{P}}}+D'\bar{Z}_{_{\hat{P}}}D)(R_{_{\hat{P}}}+D'\bar{Z}_{_{\hat{P}}}D)^{\dagger}\\
 \times(B'\bar{Z}_{_{\hat{P}}}+D'\bar{Z}_{_{\hat{P}}}C+L_{_{\hat{P}}}')\\
 =B'\bar{Z}_{_{\hat{P}}}+D'\bar{Z}_{_{\hat{P}}}C+L_{_{\hat{P}}}',
\end{array}
\right.\label{f13}
\end{eqnarray}\\
Define $P(t,T)=Z_{\hat{P}}(t,T)+\hat{P}$. It is easy to verify that $P(t,T)$ is a solution of the GDRE (\ref{f04}), is
monotonically decreasing with respect to $t$, and is bounded. Therefore, there exists a constant matrix $\bar{P}$ satisfying
\begin{eqnarray*}
\bar{P}=\lim\limits_{t\rightarrow-\infty}P(t,T)=\bar{Z}_{\hat{P}}+\hat{P}.
\end{eqnarray*}
Clearly, $\bar{P}$ is a solution of the GARE (\ref{f06}). Moreover, for  arbitrary $\hat{P}$ and $\bar{Z}_{\hat{P}}\geq 0$, it is easy to verify that $\bar{P}\geq\hat{P}$, i.e., $\bar{P}$ is the maximal solution to the GARE (\ref{f06}). This completes the proof.\\
{\bf Remark 2} \ \ The above proof implies that the solvability of the GARE (\ref{f06}) is equivalent to the solvability of the SARE (\ref{f13}).\\
To facilitate the proof of the subsequent result, we need the following remark.\\
{\bf Remark 3} \ \ If the SARE (\ref{f13}) has a solution $\bar{Z}_{\hat{P}}$, then
\begin{eqnarray}
A'\bar{Z}_{_{\hat{P}}}+\bar{Z}_{_{\hat{P}}}A+C'\bar{Z}_{_{\hat{P}}}C+Q_{_{\hat{P}}}
-M_{_{\hat{P}}}'\Omega_{_{\hat{P}}}^{\dagger}M_{_{\hat{P}}}=0,\label{f14}
\end{eqnarray}
where $M_{_{\hat{P}}}=B'\bar{Z}_{_{\hat{P}}}+D'\bar{Z}_{_{\hat{P}}}C+L_{_{\hat{P}}}',
\Omega_{_{\hat{P}}}=R_{_{\hat{P}}}+D'\bar{Z}_{_{\hat{P}}}D$.
In view of the following formula:
\begin{eqnarray*}
M_{_{\hat{P}}}'\Omega_{_{\hat{P}}}^{\dagger}M_{_{\hat{P}}}=-M_{_{\hat{P}}}'K_{_{\hat{P}}}-K_{_{\hat{P}}}'M_{_{\hat{P}}}-K'_{_{\hat{P}}}\Omega_{_{\hat{P}}}K_{_{\hat{P}}},
\end{eqnarray*}
in which $K_{_{\hat{P}}}=-\Omega_{_{\hat{P}}}^{\dagger}M_{_{\hat{P}}}$,
(\ref{f14}) can be further rewritten as follows:
\begin{eqnarray}
\bar{A}'\bar{Z}_{_{\hat{P}}}+\bar{Z}_{_{\hat{P}}}\bar{A}
+\bar{C}'\bar{Z}_{_{\hat{P}}}\bar{C}+\bar{Q}_{_{\hat{P}}}=0,\label{f15}
\end{eqnarray}
where
\begin{eqnarray}
\bar{A}=A+BK_{_{\hat{P}}},\quad \bar{C}=C+DK_{_{\hat{P}}}, \quad \nonumber\\ \bar{Q}_{_{\hat{P}}}=Q_{_{\hat{P}}}+L_{_{\hat{P}}}K_{_{\hat{P}}}
+K_{_{\hat{P}}}'L'_{_{\hat{P}}}+K_{_{\hat{P}}}'R_{_{\hat{P}}}K_{_{\hat{P}}}.\label{f16}
\end{eqnarray}
Using Extended Schur's Lemma in \cite{20Albert:69}, we have
\begin{eqnarray}
\bar{Q}_{_{\hat{P}}}&\geq& L_{_{\hat{P}}}R_{_{\hat{P}}}^{\dagger}L_{_{\hat{P}}}'+L_{_{\hat{P}}}R_{_{\hat{P}}}^{\dagger}R_{_{\hat{P}}}K_{_{\hat{P}}}+K_{_{\hat{P}}}'R_{_{\hat{P}}}R_{_{\hat{P}}}^{\dagger}L_{_{\hat{P}}}'\nonumber\\
&&+K_{_{\hat{P}}}'R_{_{\hat{P}}}R_{_{\hat{P}}}^{\dagger}R_{_{\hat{P}}}K_{_{\hat{P}}}\nonumber\\
&=&(L_{_{\hat{P}}}'+R_{_{\hat{P}}}K_{_{\hat{P}}})'R_{_{\hat{P}}}^{\dagger}
(L_{_{\hat{P}}}'+R_{_{\hat{P}}}K_{_{\hat{P}}})\geq0.\label{f17}
\end{eqnarray}\\
{\bf Definition 2} \ \ Consider the following stochastic system:
\begin{eqnarray}
\left\{
\begin{array}{lll}
dx(t)=Ax(t)dt+Cx(t)dw(t), \\
\mathcal{Y}(t)=Q^{\frac{1}{2}}x(t).
\end{array}
\right.\label{f18}
\end{eqnarray}
 $(A, C, Q^{\frac{1}{2}})$ is said to be exact detectable if, for any $T\geq0$,
 \begin{eqnarray*}
 \mathcal{Y}(t)=0,  a.s.,  \forall t\in[0,T]\quad \Rightarrow \quad\lim\limits_{t\rightarrow +\infty}E(x'(t)x(t))=0.
\end{eqnarray*}\\
%{\bf Definition 3.4} \ \ For system (\ref{f18}),
% $(A, C, Q^{\frac{1}{2}})$ is said to be exact observable, if for any $T\geq0$,
 %\begin{eqnarray*}
 %\mathcal{Y}(t)=0,  a.s.,  \forall t\in[0,T]\quad \Rightarrow \quad x_{0}=0.
%\end{eqnarray*}\\
{\bf Assumption 1 } \ \ $(A, C, Q_{\hat{P}}^{\frac{1}{2}})$ is exact detectable.\\
{\bf Lemma 2} \ \ If $(A, C, Q_{\hat{P}}^{\frac{1}{2}})$ is exact detectable, as in Assumption 1,  then $(\bar{A}, \bar{C}, \bar{Q}_{\hat{P}}^{\frac{1}{2}})$ is exact detectable, where $\bar{A}, \bar{C}, \bar{Q}_{\hat{P}}$ are defined in (\ref{f16}).\\
{\emph{Proof}:} \ \ The proof is similar to that of Proposition 1 in \cite{21Zhang:2004}, so we omit it here.\\
{\bf Theorem 2 } \ \ If $\mathcal{P}\neq\emptyset$ and Assumption 1 is satisfied, then the
closed-loop system (\ref{f01}) is mean-square stabilizable if and only if the GARE (\ref{f06}) has a solution $\bar{P}$ that is also the maximal solution to (\ref{f06}). In this case, the optimal stabilizing controller is given by
\begin{eqnarray}
u(t)=Kx(t),\label{f19}
\end{eqnarray}
where \begin{eqnarray}
K=-(R+D'\bar{P}D)^{\dagger}(B'\bar{P}+D'\bar{P}C),\label{f20}
\end{eqnarray}
 and the optimal cost functional is $J^{\ast}=E(x'_{0}\bar{P}x_{0})$.\\
{\emph{Proof}:} \ \ ``\emph{Sufficiency}": We will show that, under the conditions that $\mathcal{P}\neq\emptyset$ and Assumption 1 holds, when the GARE (\ref{f06}) has a solution $\bar{P}$, the closed-loop system (\ref{f01}) is mean-square stabilizable.\\
Let $\hat{P}\in\mathcal{P}$. From Remark 2, when the GARE (\ref{f06}) has a solution $\bar{P}$, the SARE (\ref{f13}) has a positive semidefinite solution $\bar{Z}_{\hat{P}}$, i.e., $\bar{Z}_{\hat{P}}\geq0$. Moreover, $\bar{P}=\bar{Z}_{\hat{P}}+\hat{P}$. Next we will show that system (\ref{f01}) with $u(t)=Kx(t)$ where $K$ is defined by (\ref{f20}), i.e.,
\begin{eqnarray}
dx(t)=(A+BK)x(t)dt+(C+DK)x(t)dw(t)\label{f22}
\end{eqnarray}
is mean square stabilizable. In view of the relationship
$K=K_{\hat{P}}$ which is defined in Remark 3, we can see that the stabilization for the system (\ref{f01}) with $u(t)=Kx(t)$ is equivalent to the stabilization for the system (\ref{f01}) with $u(t)=K_{\hat{P}}x(t)$.
We define the Lyapunov function candidate
as
\begin{eqnarray}
V(t,x(t))=E[x'(t)\bar{Z}_{\hat{P}}x(t)], \quad t\geq0.\label{f23}
\end{eqnarray}
In view of $\bar{Z}_{\hat{P}}\geq0$, the limits of $V(t,x(t))$ exists which can be similarly obtained from the proof of Theorem 4 in \cite{6Zhang:2017}.
 %Taking time derivative along the dynamic direction of the stochastic
%system (\ref{f01}), we have
%\begin{eqnarray}
%\dot{V}(t,x(t))&=&-E[x'(t)Q_{\hat{P}}x(t)+2x'(t)L_{\hat{P}}u(t)\nonumber\\
%&&+u'(t)R_{\hat{P}}u(t)]\nonumber\\
%&=&-E\left\{\left[
%  \begin{array}{cc}
%x(t)\\
%u(t)
%  \end{array}
%\right]'\left[
%  \begin{array}{cc}
%Q_{_{\hat{P}}}& L_{_{\hat{P}}}\\
%L_{_{\hat{P}}}' & R_{_{\hat{P}}}
%  \end{array}
%\right]\left[
%  \begin{array}{cc}
%x(t)\\
%u(t)
%  \end{array}
%\right]\right\}\\\nonumber
%&\leq&0,\label{f24}
%\end{eqnarray}
%where $Q_{\hat{P}}, R_{\hat{P}}, L_{\hat{P}}$ defined in (\ref{f07}).\\
%From (\ref{f24}), we can see $V(t,x(t))$ is nonincreasing with respect to $t$. That implies $V(t,x(t))\leq V(0,x_{0})$, i.e., $V(t,x(t))$ is bounded. Therefore, $\lim\limits_{t\rightarrow+\infty}V(t,x(t))$ exists.\\
And there exists an orthogonal matrix $U_{\hat{P}}$ such that
\begin{eqnarray}
U_{\hat{P}}'\bar{Z}_{\hat{P}}U_{\hat{P}}=\left[
  \begin{array}{cc}
0  & 0\\
0& \bar{Z}_{\hat{P}_{2}}
  \end{array}
\right], \quad \bar{Z}_{\hat{P}_{2}}>0.\label{f25}
\end{eqnarray}\\
From (\ref{f15}), we have
\begin{eqnarray}
&&U_{\hat{P}}'\bar{A}'U_{\hat{P}}U_{\hat{P}}'\bar{Z}_{\hat{P}}U_{\hat{P}}+U_{\hat{P}}'\bar{Z}_{\hat{P}}U_{\hat{P}}U_{\hat{P}}'\bar{A}U_{\hat{P}}\nonumber\\
&&+U_{\hat{P}}'\bar{C}'U_{\hat{P}}U_{\hat{P}}'\bar{Z}_{\hat{P}}U_{\hat{P}}
U_{\hat{P}}'\bar{C}U_{\hat{P}}+U_{\hat{P}}'\bar{Q}_{\hat{P}}U_{\hat{P}}=0.\label{f26}
\end{eqnarray}\\
Now we assume
\begin{eqnarray}
U_{\hat{P}}'\bar{A}U_{\hat{P}}&=&\left[
  \begin{array}{cc}
\bar{A}_{11}&\bar{A}_{12}\\
\bar{A}_{21}&\bar{A}_{22}
  \end{array}
\right],\quad U_{\hat{P}}'\bar{C}U_{\hat{P}}=\left[
  \begin{array}{cc}
\bar{C}_{11}&\bar{C}_{12}\\
\bar{C}_{21}&\bar{C}_{22}
  \end{array}
\right],\nonumber\\
U_{_{\hat{P}}}'\bar{Q}_{_{\hat{P}}}U_{_{\hat{P}}}&=&\left[
  \begin{array}{cc}
\bar{Q}_{\hat{P}_{11}}&\bar{Q}_{\hat{P}_{12}}\\
\bar{Q}_{\hat{P}_{12}}'&\bar{Q}_{\hat{P}_{22}}
  \end{array}
\right],\label{f27}
\end{eqnarray}
and by a simple calculation we can obtain
\begin{eqnarray}
  &&\begin{bmatrix}
\begin{smallmatrix}
\bar{C}_{21}'\bar{Z}_{\hat{P}_{2}}\bar{C}_{21}+\bar{Q}_{\hat{P}_{11}} & \bar{A}_{21}'\bar{Z}_{\hat{P}_{2}}+\bar{C}_{21}'\bar{Z}_{\hat{P}_{2}}\bar{C}_{22}
+\bar{Q}_{\hat{P}_{12}}\\
\bar{Z}_{\hat{P}_{2}}\bar{A}_{21}+\bar{C}_{22}'\bar{Z}_{\hat{P}_{2}}\bar{C}_{21}
+\bar{Q}_{\hat{P}_{12}} &
\bar{A}_{22}'\bar{Z}_{\hat{P}_{2}}+\bar{Z}_{\hat{P}_{2}}\bar{A}_{22}
+\bar{C}'_{22}\bar{Z}_{\hat{P}_{2}}\bar{C}_{22}+\bar{Q}_{\hat{P}_{22}}
  \end{smallmatrix}
\end{bmatrix}\nonumber\\
&&=0.\label{f28}
\end{eqnarray}\\
From Remark 3, we have $\bar{Q}_{\hat{P}}\geq0$, with $\bar{Z}_{\hat{P}_{2}}>0$, it is easy to obtain $\bar{Q}_{\hat{P}_{11}}=0$, $\bar{C}_{21}=0$. Next we will  illustrate $\bar{Q}_{\hat{P}_{12}}=0$. Actually, for any $x=U_{\hat{P}}\left[
  \begin{array}{cc}
x_{1}\\
x_{2}
  \end{array}
\right]$, where $x_{2}$ has the same dimension with $\bar{Q}_{\hat{P}_{12}}$, it has
\begin{eqnarray}
x'\bar{Q}_{\hat{P}}x&=&\left[
  \begin{array}{cc}
x_{1}\\
x_{2}
  \end{array}
\right]'U_{\hat{P}}'\bar{Q}_{\hat{P}}U_{\hat{P}}\left[
  \begin{array}{cc}
x_{1}\\
x_{2}
  \end{array}
\right]\nonumber\\
&=&x_{2}'\bar{Q}_{\hat{P}_{12}}'x_{1}+x_{1}'\bar{Q}_{\hat{P}_{12}}x_{2}
+x_{2}'\bar{Q}_{\hat{P}_{22}}x_{2}.\label{f29}
\end{eqnarray}
If $\bar{Q}_{\hat{P}_{12}}\neq0$, we can always find $x_{1}$, $x_{2}$, such that $x'\bar{Q}_{\hat{P}}x<0$ which is contrary with $\bar{Q}_{\hat{P}}\geq0$. Therefore, $\bar{Q}_{\hat{P}_{12}}=0$, $\bar{A}_{21}=0$ and $\bar{Q}_{\hat{P}_{22}}\geq0$.\\
Plugging the above results into (\ref{f28}), we have
\begin{eqnarray}
\bar{A}_{22}'\bar{Z}_{\hat{P}_{2}}+\bar{Z}_{\hat{P}_{2}}\bar{A}_{22}
+\bar{C}_{22}'\bar{Z}_{\hat{P}_{2}}\bar{C}_{22}+\bar{Q}_{\hat{P}_{22}}=0.\label{f30}
\end{eqnarray}\\
Denote $U_{\hat{P}}'x(t)=\bar{x}(t)=\left[
  \begin{array}{cc}
\bar{x}^{(1)}(t)\\
\bar{x}^{(2)}(t)
  \end{array}
\right]$, and the dimension of $\bar{x}^{(2)}(t)$ coincides with the rank of $\bar{Z}_{\hat{P}_{2}}$. Thus, (\ref{f22}) can be rewritten as
\begin{eqnarray}
U_{\hat{P}}'dx(t)=U_{\hat{P}}'\bar{A}U_{\hat{P}}U_{\hat{P}}'x(t)dt
+U_{\hat{P}}'\bar{C}U_{\hat{P}}U_{\hat{P}}'x(t)dw(t)\label{f31}
\end{eqnarray}
i.e.,
\begin{eqnarray}
d\bar{x}^{(1)}(t)&=&[\bar{A}_{11}\bar{x}^{(1)}(t)+\bar{A}_{12}\bar{x}^{(2)}(t)]dt\nonumber\\
&&+[\bar{C}_{11}\bar{x}^{(1)}(t)+\bar{C}_{12}\bar{x}^{(2)}(t)]dw(t),\label{f32}\\
d\bar{x}^{(2)}(t)&=&\bar{A}_{22}\bar{x}^{(2)}(t)dt+\bar{C}_{22}\bar{x}^{(2)}(t)dw(t).\label{f33}
\end{eqnarray}\\
Firstly, we show the stability of $(\bar{A}_{22}, \bar{C}_{22})$.\\
Applying It$\hat{o}$'s formula to $E[x'(t)\bar{Z}_{\hat{P}}x(t)]$ and taking integral from 0 to $T$, it holds that
\begin{eqnarray}
&&E[x'(T)\bar{Z}_{\hat{P}}x(T)]-E[x'_{0}\bar{Z}_{\hat{P}}x_{0}]\nonumber\\
&=&-E\int_{0}^{T}[x'(t)Q_{\hat{P}}x(t)+x'(t)L_{\hat{P}}u(t)+u'(t)L'_{\hat{P}}x(t)\nonumber\\
&&+u'(t)R_{\hat{P}}u(t)]dt\nonumber\\
&=&-E\int_{0}^{T}x'(t)\bar{Q}_{\hat{P}}x(t)dt\leq0.\label{f34}
\end{eqnarray}
Thus, it holds that
\begin{eqnarray}
&&\int_{0}^{T}E[\bar{x}^{(2)'}(t)\bar{Q}_{\hat{P}_{22}}\bar{x}^{(2)}(t)]dt
=\int_{0}^{T}E[x'(t)\bar{Q}_{\hat{P}}x(t)]dt\nonumber\\
&=&E[x'_{0}\bar{Z}_{\hat{P}}x_{0}]-E[x'(T)\bar{Z}_{\hat{P}}x(T)]\nonumber\\
&=&E[\bar{x}^{(2)'}_{0}\bar{Z}_{\hat{P}_{2}}\bar{x}^{(2)}_{0}]
-E[\bar{x}^{(2)'}(T)\bar{Z}_{\hat{P}_{2}}\bar{x}^{(2)}(T)].\label{f35}
\end{eqnarray}
Following from Lemma 3 in \cite{16Qi:2017}, we know that system $(\bar{A}_{22}, \bar{C}_{22}, \bar{Q}^{\frac{1}{2}}_{\hat{P}_{22}})$ is exact observable with $\bar{Z}_{\hat{P}_{2}}>0$.\\
Since
\begin{eqnarray}
\int_{0}^{T}E[\bar{x}^{(2)'}(t)\bar{Q}_{\hat{P}_{22}}\bar{x}^{(2)}(t)]dt
=E[\bar{x}^{(2)'}_{0}\bar{F}_{\hat{P}}(0,T)\bar{x}^{(2)}_{0}],\label{f36}
\end{eqnarray}
 where $\bar{F}_{\hat{P}}(t,T)$ satisfies the following differential equation:
\begin{eqnarray}
-\dot{\bar{F}}_{\hat{P}}(t,T)&=&\bar{Q}_{\hat{P}_{22}}+\bar{A}_{22}'\bar{F}_{\hat{P}}(t,T)
+\bar{F}_{\hat{P}}(t,T)\bar{A}_{22}\nonumber\\
&&+\bar{C}_{22}'\bar{F}_{\hat{P}}(t,T)\bar{C}_{22},\label{f37}
\end{eqnarray}
with final condition $\bar{F}_{\hat{P}}(T,T)=0$. From $\bar{Q}_{\hat{P}_{22}}\geq0$, we know $\bar{F}_{\hat{P}}(t,T)\geq0$ for $t\in[0,T]$.\\
Now we will show $\bar{F}_{\hat{P}}(0,T)>0$. If not, there exists a nonzero $y$ such that $E[y'\bar{F}_{\hat{P}}(0,T)y]=0$. Then we choose the initial state be $y$, (\ref{f36}) can be reduced to
\begin{eqnarray}
\int_{0}^{T}E[\bar{x}^{(2)'}(t)\bar{Q}_{\hat{P}_{22}}\bar{x}^{(2)}(t)]dt
=E[y'\bar{F}_{\hat{P}}(0,T)y]=0,\label{f38}
\end{eqnarray}
which implies that $\bar{Q}^{\frac{1}{2}}_{\hat{P}_{22}}\bar{x}^{(2)}(t)=0$, for any $t\in[0,T]$. With the exact observability  of system $(\bar{A}_{22}, \bar{C}_{22}, \bar{Q}^{\frac{1}{2}}_{\hat{P}_{22}})$, we can obtain that $y=0$, which is contrary with $y\neq0$. Therefore, we have $\bar{F}_{\hat{P}}(0,T)>0$.\\
Via a time shift of $t$, combining (\ref{f36}), we have
\begin{eqnarray}
&&\int_{t}^{t+T}E[\bar{x}^{(2)'}(s)\bar{Q}_{\hat{P}_{22}}\bar{x}^{(2)}(s)]ds\nonumber\\
&=&E[\bar{x}^{(2)'}(t)\bar{F}_{\hat{P}}(0,T)\bar{x}^{(2)}(t)]\nonumber\\
&=&V_{2}(t,\bar{x}^{(2)}(t))-V_{2}(t+T,\bar{x}^{(2)}(t+T)),\label{f39}
\end{eqnarray}
where $V_{2}(t,x(t))=E[x'(t)\bar{Z}_{\hat{P}_{2}}x(t)]$ and similar to the proof of the convergence of $V(t,x(t))$, we know $\lim\limits_{t\rightarrow+\infty}V_{2}(t,x(t))$ exists. Taking limitation on both sides of (\ref{f39}), we have that
\begin{eqnarray} \lim\limits_{t\rightarrow+\infty}E[\bar{x}'^{(2)}(t)\bar{x}^{(2)}(t)]=0 \label{f40}
 \end{eqnarray}
 ,i.e., system $(\bar{A}_{22}, \bar{C}_{22}, \bar{Q}^{\frac{1}{2}}_{\hat{P}_{22}})$ is mean square stabilizable.\\
Next we will illustrate the stability of $(\bar{A}_{11}, \bar{C}_{11})$. Let $\bar{x}^{(2)}(0)=0$, from (\ref{f33}) we know $\bar{x}^{(2)}(t)=0$, $t\geq0$. Now (\ref{f32}) can be rewritten as
\begin{eqnarray}
d\bar{x}^{(1)}(t)=\bar{A}_{11}\bar{x}^{(1)}(t)dt+\bar{C}_{11}
\bar{x}^{(1)}(t)dw(t).\label{f41}
\end{eqnarray}
Under the condition of $\bar{x}^{(2)}(t)=0$, it is easy to see that
\begin{eqnarray}
E[x'(t)\bar{Q}_{\hat{P}}x(t)]=E[\bar{x}'^{(2)}(t)\bar{Q}_{\hat{P}_{22}}
\bar{x}^{(2)}(t)]=0.\label{f42}
\end{eqnarray}
Hence, from the exact detectability of $(\bar{A}, \bar{C}, \bar{Q}^{\frac{1}{2}}_{\hat{P}})$ and $\bar{x}^{(2)}(t)=0$ we know that
\begin{eqnarray}
&&\lim\limits_{t\rightarrow+\infty}E[\bar{x}'^{(1)}(t)\bar{x}^{(1)}(t)]\nonumber\\
&=&\lim\limits_{t\rightarrow+\infty}{E[\bar{x}'^{(1)}(t)\bar{x}^{(1)}(t)]
+E[\bar{x}'^{(2)}(t)\bar{x}^{(2)}(t)]}\nonumber\\
&=&\lim\limits_{t\rightarrow+\infty}E[\bar{x}'(t)\bar{x}(t)]
=\lim\limits_{t\rightarrow+\infty}E[x'(t)x(t)]=0,\label{f43}
\end{eqnarray}
which means $(\bar{A}_{11}, \bar{C}_{11})$ is mean square stable.\\
Secondly, we will show system (\ref{f01}) with controller $u(t)=Kx(t)$ is stabilizable in the mean square sense. In fact, we denote $\bar{\mathcal{A}}=\left[
 \begin{array}{cc}
\bar{A}_{11} &0 \\
0&\bar{A}_{22}
  \end{array}
\right], \bar{\mathcal{C}}=\left[
  \begin{array}{cc}
\bar{C}_{11} &0 \\
0&\bar{C}_{22}
  \end{array}
\right]$. Thus (\ref{f32}) and (\ref{f33}) can be rewritten as below
\begin{eqnarray}
d\bar{x}(t)&=&\big\{\bar{\mathcal{A}}\bar{x}(t)+\left[
  \begin{array}{cc}
\bar{A}_{12}\\
0
  \end{array}
\right]\bar{u}(t)\big\}dt\nonumber\\
&&+\big\{\bar{\mathcal{C}}\bar{x}(t)+\left[
  \begin{array}{cc}
\bar{C}_{12}\\
0
  \end{array}
\right]\bar{u}(t)\big\}dw(t),\label{f44}
\end{eqnarray}
where $\bar{u}(t)$ is the solution to equation (\ref{f33}) with initial condition $\bar{u}(0)=\bar{x}^{(2)}(0)$. The stability of $(\bar{A}_{11}, \bar{C}_{11})$ and $(\bar{A}_{22}, \bar{C}_{22})$ as shown above indicates that $(\bar{\mathcal{A}}, \bar{\mathcal{C}})$ is stable in the mean square sense. From (\ref{f40}) it is easy to see that $\lim\limits_{t\rightarrow+\infty}E[\bar{u}'(t)\bar{u}(t)]=0$ and $\int_{0}^{\infty}E[\bar{u}'(t)\bar{u}(t)]dt<+\infty$. Applying the corresponding results in \cite{22Hinrichsen:98}, we have that there exists a constant $c_{0}$ satisfying
\begin{eqnarray}
\int_{0}^{\infty}E[\bar{x}'(t)\bar{x}(t)]dt<c_{0}\int_{0}^{\infty}
E[\bar{u}'(t)\bar{u}(t)]dt<+\infty.\label{f45}
\end{eqnarray}
Hence, $\lim\limits_{t\rightarrow+\infty}E[\bar{x}'(t)\bar{x}(t)]=0$ can be verified from (\ref{f45}) and
\begin{eqnarray}
\lim\limits_{t\rightarrow+\infty}E[x'(t)x(t)]
=\lim\limits_{t\rightarrow+\infty}E[\bar{x}'(t)\bar{x}(t)]=0.\label{f46}
\end{eqnarray}
In conclusion, system (\ref{f01}) can be stabilizable with controller $u(t)=Kx(t)$ in the mean square sense.\\
Finally, we will give the optimal controller and the associated cost functional. Applying It$\hat{o}$'s
formula to $x'(t)\bar{P}x(t)$ with the GARE (\ref{f06}), we have
\begin{eqnarray*}
2J_{T}&=&E\int_{0}^{T}\bigg\{x'(t)[A'\bar{P}+\bar{P}A+C'\bar{P}C]x(t)+u'(t)[B'\bar{P}\\
&&+D'\bar{P}C]x(t)+x'(t)[\bar{P}B+C'\bar{P}D]u(t)\\
&&+u'(t)D'\bar{P}Du(t)\bigg\}dt+E(x'_{0}\bar{P}x_{0})-E[x'(T)\bar{P}x(T)]\\
&=&E\int_{0}^{T}[u(t)-Kx(t)]'(R+D'\bar{P}D)[u(t)-Kx(t)]dt\\
&&+E(x'_{0}\bar{P}x_{0})-E[x'(T)\bar{P}x(T)].
\end{eqnarray*}
With the mean square stabilizable of system (\ref{f01}), let $T\rightarrow+\infty$ in the above equation, the infinite cost functional can be reformulated as
\begin{eqnarray}
J&=&E(x'_{0}\bar{P}x_{0})+E\int_{0}^{\infty}[u(t)-Kx(t)]'(R+D'\bar{P}D)\nonumber\\
&&\cdot[u(t)-Kx(t)]dt.\label{f47}
\end{eqnarray}
Thus the optimal controller is $u(t)=Kx(t)$ where $K$ as (\ref{f20}). In addition, the optimal cost functional is $J^{\ast}=E(x'_{0}\bar{P}x_{0})$.\\
``\emph{Necessity}": On the other hand, we should illustrate that if $\mathcal{P}\neq\emptyset$, when the closed-loop system (\ref{f01}) is mean square stabilizable, then the GARE (\ref{f06}) has a solution $\bar{P}$, which is also the maximal solution to the GARE (\ref{f06}).
In fact, the existence of solutions to the GARE (\ref{f06}) have been obtained in Theorem 1.\\
The proof is complete.\\
{\bf Remark 4} \ \  The key technique of sufficient proof is that in terms of the equivalence between the GARE and the SARE, the stabilization of the general indefinite case is reduced to the definite one whose stabilization is studied by Lyapunov functional defined with the optimal cost functional.
\section{An example}
In this section, we give an example to illustrate the result obtained by Theorem 2.\\
Consider two-dimensional system (\ref{f01}) and cost functional (\ref{f02}) with the following parameters
$A=\left[
  \begin{array}{cc}
0.01 &0\\
0&-0.1
  \end{array}
\right]
, B=\left[
  \begin{array}{cc}
0.2\\
0
  \end{array}
\right]
, C=\left[
  \begin{array}{cc}
-0.1&0\\
0&0.1
  \end{array}
\right], D=\left[
  \begin{array}{cc}
0.6\\
0
  \end{array}
\right], Q=\left[
  \begin{array}{cc}
0.5&0\\
0&-1
  \end{array}
\right], R=-0.05$ with initial state $x_{0}=\left[
  \begin{array}{cc}
-0.01\\
0.1
  \end{array}
\right]
$, noting that $Q, R$ are indefinite matrices.  By solving the LMI in set $\mathcal{P}$ with Matlab tool, the solution to GARE (\ref{f06}) with the given parameters can be obtained as $\bar{P}=\left[
  \begin{array}{cc}
20.143&0\\
0&-5.2632
  \end{array}
\right]
$.  Therefore,  the gain matrix $K$ in (\ref{f20}) can be easily obtained as $K=\left[
  \begin{array}{cc}
-0.2891&0
  \end{array}
\right]
$. Furthermore,  Assumption 1 is also satisfied due to $Q_{\bar{P}}=\left[
  \begin{array}{cc}
0.9230 &0\\
0&0.1053
  \end{array}
\right]
$.
Accordingly, by Theorem 2 we know that the system (\ref{f01}) with the above parameters is mean square stabilizable. The simulation result is shown in Fig. 1. It can be seen that the state $x(t)$ (Fig. 3 ) is stabilized with the optimal controller (Fig. 2), as expected.
  \begin{figure}[htbp]
  \begin{center}
  \includegraphics[width=0.42\textwidth]{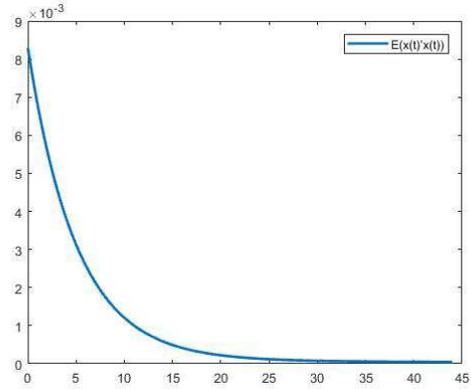}
  \caption{Simulations for the state trajectory $E[x'(t)x(t)]$.} \label{fig:digit}
  \end{center}
\end{figure}
\begin{figure}[htbp]
  \begin{center}
  \includegraphics[width=0.42\textwidth]{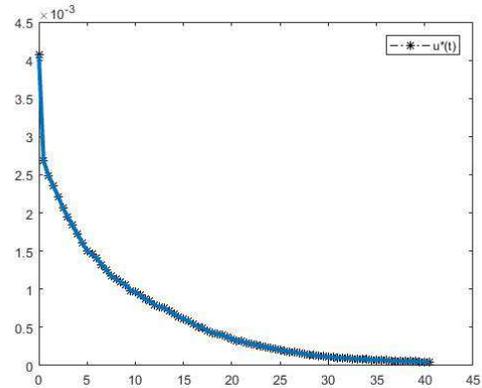}
  \caption{Optimal control.} \label{fig:digit}
  \end{center}
\end{figure}
\begin{figure}[htbp]
  \begin{center}
  \includegraphics[width=0.42\textwidth]{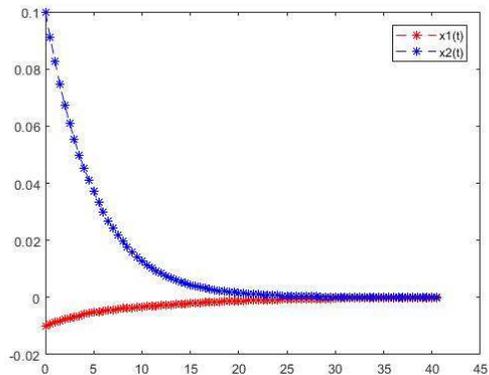}
  \caption{Optimal states.} \label{fig:digit}
  \end{center}
\end{figure}
\section{Conclusions}         % Sections and subsections are supported
In this paper, we mainly discussed the stabilization problem for It$\hat{o}$ stochastic system, whose control and state weighting matrices in the cost functional are indefinite.  The convergence of the GDRE which involves a matric pseudo-inverse and two additional equality/inequality constraints was studied.  And the infinite horizon optimal controller was  obtained accordingly. Finally, in terms of the equivalence between the GARE and the SARE,  the stabilization of the general indefinite case was reduced to the definite one whose stabilization is studied by Lyapunov functional defined with the optimal cost functional  subject to the SARE. The contents in this paper are an extension and improvement of the previous works \cite{19Rami:2001}.

\bibliographystyle{plain}        % Include this if you use bibtex
\bibliography{autosam}           % and a bib file to produce the
                                 % bibliography (preferred). The
                                 % correct style is generated by
                                 % Elsevier at the time of printing.

%\appendix
%\section{A summary of Latin grammar}    % Each appendix must have a short title.
%\section{Some Latin vocabulary}         % Sections and subsections are supported
                                        % in the appendices.
\end{document}